\documentclass[12 pt,a4paper]{article}
\usepackage{cite}
\usepackage{amsmath,dsfont,amsfonts,latexsym,amscd,amssymb}
\usepackage{extsizes}
\usepackage{graphicx}
\usepackage{float}
\usepackage[linkcolor=blue, urlcolor=blue, citecolor=blue, colorlinks, bookmarks]{hyperref}
\usepackage{amsfonts}
\usepackage[left=2cm,right=2cm,top=2cm,bottom=2cm]{geometry}
\usepackage{setspace}
\newtheorem{Theorem}{Theorem}[section]

\newtheorem{Remark}{Remark}

\newtheorem{Lemma}{Lemma} 
\numberwithin{Theorem}{section} \numberwithin{equation}{section}

\begin{document}
	\title{\textbf{On Hypersurface of a Finsler space subjected to \textsl{h-}Matsumoto change}}
	\author{M.$\,$K.$\,$\textsc{Gupta} and Suman \textsc{Sharma}\\
		\normalsize{Department of Mathematics}\\[-3mm]
		\normalsize{Guru Ghasidas Viswavidyalaya, Bilaspur (C.G.), India}\\[-3mm]
		\small{E-mail: mkgiaps@gmail.com, sharma.suman209@gmail.com}}
	\date{}
	\maketitle{}
	\begin{abstract} Recently, we have studied the Finsler space with \textsl{h-}Matsumoto change and found Cartan connection for the transformed space \cite{first paper}. In this paper, we have discussed certain geometrical properties of the hypersurface of a Finsler space for the \textsl{h-}Matsumoto change.
	\end{abstract}
	\section{Introduction}In 1984, C. Shibata \cite{shibata1984invariant} introduced a change of Finsler metric called $\beta$-change. An important class of $\beta$-change is Matsumoto change which is given by 
	$\overline{L}(x,y) = \frac{L^{2}}{L-\beta}$,\, where \,$\beta(x,y)= b^{}_{i}(x)y^{i}$ is one form on $M^{n}$. The concept of an \textsl{h-}vector $b^{}_{i}(x,y)$ has been introduced by H. Izumi \cite{izumi1980conformal}, which is v-covariant constant and satisfies $LC^{h}_{ij} b_{h}= \rho h_{ij}$, where $\rho$ is a non-zero scalar function. He showed that the scalar $\rho$ is independent of directional arguments. M.$\,$K.$\,$Gupta and P.$\,$N.$\,$Pandey \cite{gupta2015finsler} proved that the scalar $\rho$ is constant	if the \textsl{h-}vector $b_{i}$ is gradient.
B. N. Prasad \cite{prasad} discussed the Cartan connection of Finsler space whose metric is given by \textsl{h}-Randers change of a Finsler metric. In 2016, M. K. Gupta and A. K. Gupta \cite{gupta2016hypersurface} studied Finsler space subjected to \textsl{h-}exponential change.  Recently, we have obtained \cite{first paper} the relation between the Cartan connection of Finsler space $F^{n}=(M^{n},L)$ and $\overline{F}^{n}= (M^{n},\overline{L})$, where $\overline{L}(x,y)$ is obtained by the transformation \\[-8mm]
\begin{equation}\label{eq3}
	\overline{L}(x,y)= \frac{L^{2}(x,y)}{L(x,y)-b_{i}(x,y)y^{i}},
\end{equation} \\[-8mm]
and $b_{i}$ is an $\textsl{h-}$vector in $(M^{n},L)$.\\
  \noindent ``A Hypersurface is a generalization of the concept of hyperplane''. The theory of hypersurface in a Finsler space  has been first considered by  E.$\,$Cartan \cite{Cartan} from two points of view, \textit{i.e.}\\
  \newpage
 \begin{enumerate}
 	\item[(i)] A hypersurface as the whole of tangent line-elements and then it is also a Finsler space.\\[-9mm]
 	\item[(ii)] A hypersurface as the whole of normal line-elements and then it is a Riemannian space.\\[-9mm]
 \end{enumerate}
 Three kinds of hyperplane has been introduced by A Rapcs$\grave{a}$k \cite{Rapcsak} while M. Matsumoto \cite{Intrinsic} has classified them and developed a systematic theory of Finslerian hypersurfaces. M. K. Gupta and P. N. Pandey \cite{Gupta h2018} discussed the hypersurfaces of a Finsler space whose metric is given by Kropina change with an \textsl{h-}vector. The hypersurface of a Finsler space subjected to an \textsl{h-}exponential change of metric has been studied by Gupta and Gupta \cite{gupta2017h}.\\
In the present paper, we discuss the geometrical properties of hypersurface of a Finsler space subjected to the \textsl{h-}Matsumoto change given by \eqref{eq3}.\\
The terminologies and notations are referred to Matsumoto \cite{Matsumoto Book}.
\section{Preliminaries} Let $F^{n}=(M^{n},L)$ be an $n$-dimensional Finsler space equipped with the fundamental function $L(x,y)$, satisfying the requisite conditions \cite{Matsumoto Book}. The normalized supporting element, angular metric tensor, metric tensor  and Cartan tensor are defined by $l_{i}= \dot{\partial}_{i}L$,$~$ $h_{ij}= L\dot{\partial_{i}}\dot{\partial_{j}}L$,\, $g_{ij}=\frac{1}{2}\dot{\partial_{i}}\dot{\partial_{j}}L^{2}$ and $C_{ijk} = \frac{1}{2}\dot{\partial}_{k} g_{ij}$ respectively. The Cartan connection in $F^{n}$ is given as $C\Gamma= (F^{i}_{jk}, G^{i}_{j}, C^{i}_{jk})$\vspace{.3cm}.\\[-1mm]
A hypersurface $M^{n-1}$ of the underlying smooth manifold $M^{n}$ may be parametrically represented by the equation $x^{i}= x^{i}(u^{\alpha})$, where $u^{\alpha}$ are Gaussian coordinates on $M^{n-1}$ and Greek indices run from $1$ to $n-1$. We assume that the matrix of projection factors $B^{i}_{\alpha}= \frac{\partial{x}^{i}}{\partial{u}^{\alpha}}$ is of rank $n-1$. If the supporting element $y^{i}$ at a point $u=(u^{\alpha})$ of $M^{n-1}$ is assumed to be tangent to $M^{n-1}$, we may then write $y^{i}=B^{i}_{\alpha}(u)v^{\alpha}$, so that $v=(v^{\alpha})$ is thought of as the supporting element of $M^{n-1}$ at the point $u^{\alpha}$. Since the function $\underline{L}(u,v)= L\Big(x(u),y(u,v)\Big)$ gives rise to a Finsler metric on $M^{n-1}$, we get an $(n-1)$-dimensional Finsler space $F^{n-1}=(M^{n-1},\underline{L}(u,v))$.\\ 
A unit normal vector $B^{i}(u,v)$  at each point $u^{\alpha}$ of $F^{n-1}$ is defined by \cite{Intrinsic},\\[-6mm]
\begin{equation}\label{eq6}
	g_{ij}B^{i}_{\alpha}B^{j}= 0, \hspace{1.2cm} g_{ij}B^{i}B^{j}=1.
\end{equation}
 If the inverse projection  of $B^{i}_{\alpha}$ is $B^{\alpha}_{i}(u,v)$, then we have\\[-6mm]
\begin{equation}\label{eq7}
	B^{\alpha}_{i}= g^{\alpha\beta}g_{ij}B^{j}_{\beta},
\end{equation}
where $g^{\alpha\beta}$ is the inverse of the metric tensor $g_{\alpha\beta}$ of $F^{n-1}$.\\ 
In view of equation \eqref{eq6} and \eqref{eq7}, we get\\[-4mm]
\begin{equation}\label{eq8}
	B^{i}_{\alpha} B^{\beta}_{i }= \delta^{\beta}_{\alpha},\hspace{.6cm} B^{i}_{\alpha}B_{i}=0,\hspace{.6cm} B^{i}B^{\alpha}_{i}=0,\hspace{.6cm} B^{i}B_{i}=1,
\end{equation}\\[-8mm]
and further\\[-12mm] 
\begin{equation}\label{3}
	B^{i}_{\alpha}B^{\alpha}_{j}+B^{i}B_{j}= \delta^{i}_{j}.
\end{equation}
The second fundamental \textsl{h-}vector $H_{\alpha\beta}$ and the normal curvature vector $H_{\alpha}$ for the induced Cartan connection IC$\Gamma$=$(F^{\alpha}_{\beta\gamma},G^{\alpha}_{\beta},C^{\alpha}_{\beta\gamma})$ on $F^{n-1}$ are given by   \\[-3mm]
\begin{equation}\label{eq9}
	H_{\alpha\beta}=B_{i}(B^{i}_{\alpha\beta}+F^{i}_{jk}B^{j}_{\alpha}B^{k}_{\beta})+ M_{\alpha}H_{\beta},
\end{equation}
and\\[-12mm]
\begin{equation}\label{eq10}
	H_{\alpha}=B_{i}(B^{i}_{0\alpha}+G^{i}_{j}B^{j}_{\alpha}),
\end{equation}
where\\[-8mm]
\begin{equation*}
	M_{\alpha}=C_{ijk}B^{i}_{\alpha}B^{j}B^{k},\hspace{.7cm} B^{i}_{\alpha\beta}=\frac{\partial^{2}x^{i}}{\partial u^{\alpha}\partial u^{\beta}},\hspace{.7cm}  B^{i}_{0\alpha}=B^{i}_{\beta \alpha}v^{\beta}.
\end{equation*}
The equation \eqref{eq9} and \eqref{eq10} gives,\\[-6mm]
 \begin{equation}\label{eq11}
 	H_{0\alpha}=H_{\beta\alpha}v^{\beta}=H_{\alpha} \quad \text{and} \hspace{.7cm}
 	H_{\alpha0}=H_{\alpha\beta}v^{\beta}=H_{\alpha}+M_{\alpha}H_{0}.
 \end{equation}
The second fundamental \textsl{v-}tensor $M_{\alpha\beta}$ is defined as\\[-6mm]
\begin{equation}\label{eq12}
	M_{\alpha\beta}=C_{ijk}B^{j}_{\alpha}B^{k}_{\beta}B^{i}.
\end{equation}
The relative \textsl{h-} and \textsl{v-}covariant derivative of $B^{i}_{\alpha}$ and $B^{i}$ are given by\\[-6mm]
 \begin{equation}\label{eq13}
 	B^{i}_{\alpha|\beta}=H_{\alpha\beta}B^{i}, \hspace{.5cm} B^{i}_{\alpha}|^{}_{\beta}=M_{\alpha\beta}B^{i},\hspace{.5cm} B^{i}\,_{^{}|\beta}= -H_{\alpha\beta}B^{\alpha}_{j}g^{ij},\hspace{.5cm} B^{i}|_{\beta}= -M_{\alpha\beta}B^{\alpha}_{j}g^{ij}.
\end{equation}
The relative \textsl{h-} and \textsl{v-}covariant derivatives of a covarient vector field $X_{i}$ are given by,\\[-6mm]
\begin{equation}\label{eq14}
	X_{i|\beta}= X_{i|j}B^{j}_{\beta}+ X^{}_{i}|_{j}B^{j}H_{\beta}, \hspace{.5cm} \qquad X^{}_{i}|_{\beta}=X^{}_{i}|_{j}B^{j}_{\beta}.
\end{equation}
 Different kinds of hyperplane and their characteristic conditions are classified by Matsumoto \cite{Intrinsic} which are given in the following Lemmas\\[-9mm]
\begin{Lemma}\label{L1} 
A hypersurface $F^{n-1}$ is a hyperplane of the first kind if and only if $H_{\alpha}=0$ or equivalently $H_{0}=0$.\\[-10mm]
\end{Lemma} 
\begin{Lemma}\label{L2}
	A hypersurface $F^{n-1}$ is a hyperplane of the second kind if and only if $H_{\alpha\beta}=0$.\\[-10mm]
\end{Lemma} 
\begin{Lemma}\label{L3}
	A hypersurface $F^{n-1}$ is a hyperplane of the third kind if and only if $H_{\alpha\beta} = 0 = M_{\alpha\beta}$.
\end{Lemma}
\section{The Finsler space $\overline{F}^{n}=(M^{n},\overline{L})$}
Let $F^{n}=(M^{n},L)$ be the Finsler space equipped with the fundamental function $L(x,y)$ given by \eqref{eq3}, where $\beta=b_{i}(x, y)y^{i}$, $b_{i}$ is \textsl{h-}vector defined by\\[-6mm]
\begin{equation}\label{eq1}
	{(i)} \hspace{.2cm}b_{i}|_{k}=0, \hspace{1.6cm} {(ii)}\hspace{.2cm}LC^{h}_{ij} b_{h}= \rho h_{ij}, \quad \rho \neq0\,.
\end{equation}
 Thus from the above definition, we have\\[-6mm]
\begin{equation}\label{eq2}
	L\dot{\partial}_{j}b_{i}= \rho h_{ij}.
\end{equation}
 In this paper, the geometric objects corresponding to $\overline{F}^{n}$ is denoted by bar over the quantity.
 \newpage
\noindent We have obtained the normalized supporting element and the angular metric tensor of $\overline{F}^{n}$ as \cite{first paper}
\begin{equation}\label{eq15}
	\overline{l}_{i}= \frac{\tau}{(\tau-1)}l_{i}+ \frac{\tau^{2}}{(\tau-1)^{2}}m_{i},
\end{equation}
and
\begin{equation}\label{eq16}
	\overline{h}_{ij}= \frac{\tau^{2}(\tau+\rho\tau-2)}{(\tau-1)^{3}}h_{ij}+ \frac{2\tau^{4}}{(\tau-1)^{4}}m_{i}m_{j},
\end{equation}
where $\tau=\frac{L}{\beta}$,\quad $m_{i}=b_{i}-\frac{1}{\tau}l_{i}$.\\
The metric tensor and Cartan tensor of the transformed space are derived as follows \cite{first paper} \\[-3mm]
\begin{equation}\label{eq17}
	\overline{g}_{ij}= pg_{ij}+ p_{1}l_{i}l_{j}+ p_{2}(m_{i}l_{j}+m_{j}l_{i})+p_{3}m_{i}m_{j},
\end{equation}
and
\begin{equation}\label{eq18}
	\overline{C}_{ijk}= pC_{ijk}+ V_{ijk},
\end{equation}
where \begin{equation*}
	p=\frac{\tau^{2}(\tau+\rho\tau-2)}{(\tau-1)^{3}}, \hspace{.4cm} p^{}_{1}=\frac{\tau^{2}(1-\rho\tau)}{(\tau-1)^{3}},\hspace{.4cm} p^{}_{2}=\frac{\tau^{3}}{(\tau-1)^{3}},\hspace{.4cm} p^{}_{3}= \frac{3\tau^4}{(\tau-1)^4},
\end{equation*}
\begin{equation*}\label{eq19}
	 \,\,V_{ijk}= K_{1}(h_{ij}m_{k}+h_{jk}m_{i}+h_{ki}m_{j})+K_{2}m_{i}m_{j}m_{k},
\end{equation*}
\begin{equation*}
	K_{1}= \frac{\tau^{3}(\tau+3\rho\tau-4)}{2L(\tau-1)^4 },\hspace{.3cm} K_{2}= \frac{6\tau^{4}}{\beta(\tau-1)^5}.
\end{equation*}\\[-12mm]
\begin{Remark}\label{R3.1}
	$V_{ijk}$ is an indicatory tensor and satisfies;
	\begin{enumerate}
		\item[\emph{(i)}] $V_{ijk}m^{i}=(2K_{1}+K_{2}m^{2})m_{j}m_{k}+K_{1}m^{2}h_{jk}$,
		\item [\emph{(ii)}] $V_{ijk}g^{ir}=K_{1}\left( h_{jk}m^{r}+h^{r}_{k}m_{j}+h^{r}_{j}m_{k}\right)+K_{2}m_{j}m_{k}m^{r}$,
	\end{enumerate}
where $h^{r}_{k}= h_{ik}g^{ir}$.
\end{Remark}
The inverse metric tensor of $\overline{F}^{n}$ is obtained as follows \cite{first paper}
\begin{equation}\label{eq20}
	\overline{g}^{ij}=qg^{ij}+ q^{}_{1}l^{i}l^{j}+ q^{}_{2}(l^{i}m^{j}+l^{j}m^{i})+ q^{}_{3}m^{i}m^{j},
\end{equation}\\[-12mm]
where\\[-6mm]
\begin{equation*}
	\quad q= \frac{1}{p}\,,\qquad q_{1}= \frac{-1}{2}\Big[\frac{p_{1}{p^{}_{3}}-p^{2}_{2}}{(p_{1}+p)p_{3}-p_{2}^{2}}+\frac{2p^{2}p^{2}_{2}p_{3}}{(3p+2p_{3}m^{2})\{(p_{1}+p)p_{3}-p_{2}^{2}\}^{2}}\Big] ,
\end{equation*}	
\begin{equation*}
	q_{2}=\frac{-2p_{2}p_{3}}{(3p+2p_{3}m^{2})\{(p_{1}+p)p_{3}-p_{2}^{2}\}}, \qquad q_{3}=\frac{-2p_{3}}{p(3p+2p_{3}m^{2})}\,,
\end{equation*}	
and $m$ is the magnitude of the vector $m^{i}=g^{ij}m_{j}$.\\
We have obtained the relation between the Cartan connection coefficients $F^{i}_{jk}$ and $\overline{F}^{i}_{jk}$ as \cite{first paper}\\[-3mm]
\begin{equation}\label{eq21}
	\overline{F}^{i}_{jk}= F^{i}_{jk}+ D^{i}_{jk}.
\end{equation}\\[-8mm]
The difference tensor $D^{i}_{jk}$ is given by \\[-1mm]
\begin{equation}\begin{split}\label{eq26}
		{D}^{\,i}_{jk}=\overline{g}^{\,is}&\Big\{Q_{j}F_{sk}+Q_{s}E_{kj}+Q_{k}F_{js}+p\left(C_{jkm}D^{m}_{s}-C_{skm}D^{m}_{j}-C_{jsm}D^{m}_{k}\right)\\& \qquad+V_{jkm}D^{m}_{s}-V_{skm}D^{m}_{j}-V_{jsm}D^{m}_{k}+B_{js}\beta_{k}-B_{jk}\beta_{s}+B_{sk}\beta_{j}\\& \qquad+\frac{p_{2}}{2}\left(\rho_{k} h_{js}-\rho_{s} h_{jk}+\rho_{j} h_{sk}\right) \Big\},
\end{split}\end{equation}
where \\[-5mm]
\begin{equation}\label{eq25}
		\quad {D}^{\,i}_{j}= \overline{g}^{\,ir}\Big\{-2D^{m}(pC_{mrj} + V_{mrj}) + Q_{r}E_{j0} + E_{00}B_{rj}+p_{2}LF_{rj} + Q_{j}F_{r0} + \frac{1}{2}p_{2}\rho^{}_{k}h_{rj}y^{k}\Big\},
\end{equation}	
\begin{equation}\label{eq24}
	\quad {D}^{\,i}=\frac{1}{2}\overline{g}^{\,is}\big\{Q_{s}E_{00} + 2p_{2}LF_{s0}\big\}\,,\\
\end{equation}
and \\[-6mm]
\begin{equation*}
\quad	Q_{r}=(p_{2}l_{r}+p_{3}m_{r}), \quad B_{rj}= K_{1}h_{rj}+K_{2}m_{r}m_{j}
\end{equation*}
\begin{equation}\label{eq27}
	2E_{ij}= b_{i|j}+b_{j|i}, \quad 2F_{ij}=b_{i|j}-b_{j|i}
\end{equation}
	\begin{equation*}
		\beta_{j}=\beta_{|j}, \quad \rho^{}_{k}=\rho^{}_{|k}=\partial_{k}\rho.
	\end{equation*}
  The zero `0' in subscript is denoted for the contraction by $y^{i}$, for example, $F_{ij}y^{j}=F_{i0}$\,.\\
\noindent If the \textsl{h-}vector $b_{i}$ is parallel, \textit{i.e.} $b_{i|j}=0$, then the Cartan connection coefficient for both spaces are equivalent. Moreover then the Berwald connection coefficient for both the spaces are also identical\cite{first paper}.\\[-6mm]
\section{Hypersurface $\overline{F}^{n-1}$ of the space $\overline{F}^n$}
Let $F^{n-1}=(M^{n-1},\underline{L}(u,v))$ be a Finslerian hypersurface of the space $F^{n}$. The functions $B^{i}_{\alpha}(u)$ may be considered as the components of $(n-1)$ linearly independent vectors tangent to $F^{n-1}$. Let $B^{i}$ be the unit normal vector at a point of $F^{n-1}$. Then the unit normal vector $\overline{B}^{i}(u,v)$ of $\overline{F}^{n-1}$ is uniquely determined by \\[-8mm]
\begin{equation}\label{eq28}
\overline{g}_{ij}B^{i}_{\alpha}\overline{B}^{j}=0,\quad \overline{g}_{ij}\overline{B}^{i}\overline{B}^{j}=1.	 
\end{equation}
The inverse projection factors $\overline{B}^{\alpha}_{i}$ are uniquely defined along $\overline{F}^{n-1}$ by \\[-4mm]
\begin{equation}\label{eq29}
 \overline{B}^{\alpha}_{i}= \overline{g}^{\alpha\beta}g_{ij}B^{j}_{\beta},
\end{equation} \\[-8mm]
where $\overline{g}^{\alpha\beta}$ is the inverse of the metric tensor $\overline{g}_{\alpha\beta}$	of $\overline{F}^{n-1}$.\\From \eqref{eq29}, it follow that\\[-6mm]
\begin{equation}\label{eq30}
	B^{i}_{\alpha}\overline{B}^{\beta}_{i}= \delta^{\beta}_{\alpha},\quad B^{i}_{\alpha}\overline{B}_{i}=0, \quad \overline{B}^{i}\overline{B}^{\alpha}_{i}=0,\quad \overline{B}^{i}\overline{B}_{i}=1.
\end{equation}
   Transvecting equation \eqref{eq6} by $v^{\alpha}$ and using $B^{i}_{\alpha}v^{\alpha}=y^{i}$, we get\\[-6mm]
	\begin{equation}\label{eq31}
		y_{j} B^{j}=0.
	\end{equation}
 Equation \eqref{eq17} is contracting by $B^{i}B^{j}$ and using \eqref{eq28} and \eqref{eq31} we have,\\[-6mm]
\begin{equation}\label{eq32}
	\overline{g}_{ij}B^{i}B^{j}= p+ p_{3}(m_{i}B^{i})^{2},
\end{equation}\\[-6mm]
which shows that ${B^{i}}/{\sqrt{p + p_{3}(m_{i}B^{i})^{2}}}$ is a unit normal vector. Again contracting \eqref{eq17} by $B^{i}_{\alpha}B^{j}$ and using 
\eqref{eq6}, \eqref{eq31}, we obtain\\[-6mm]
\begin{equation}\label{eq33}
	\overline{g}_{ij}B^{i}_{\alpha}B^{j}= (p_{2}l_{i}+ p_{3}m_{i})B^{i}_{\alpha}\,(B^{j}m_{j}). 
\end{equation}
The above equation shows that the vector $B^{j}$ is normal to $\overline{F}^{n-1}$ if and only if \\[-6mm]
\begin{equation*}
	(p_{2}l_{i}+ p_{3}m_{i})B^{i}_{\alpha}(B^{j}m_{j})=0\,.
\end{equation*}
This implies atleast one of the following holds\\[-6mm]
\begin{equation*}
	{(i)} \hspace{.3cm}	(p_{2}l_{i}+ p_{3}m_{i})B^{i}_{\alpha}=0  \quad \qquad {(ii)} \hspace{.3cm}B^{j}m_{j}=0\,.
\end{equation*} 
{(i)} on transvecting by $v^{\alpha}$ gives $L=0$, which is not possible. Therefore {(ii)} holds, \textit{i.e.}\\[-6mm]
\begin{equation}
	B^{j}m_{j}=0\,,
\end{equation}\\[-8mm]
which, in view of \eqref{eq31}, can be equivalently written as \\[-6mm]
\begin{equation}\label{eq34}
	B^{j}b_{j}=0.
\end{equation}
This shows that the vector $B^{j}$ is normal to $\overline{F}^{n-1}$ if and only if $b_{j}$ is tangent to $\overline{F}^{n-1}$. In view of equation \eqref{eq32}, \eqref{eq33} and \eqref{eq34}  we can say that ${B^{i}}/{\sqrt{p}}$ is a unit normal vector of $\overline{F}^{n-1}$\, \textit{i.e.}
\begin{equation}\label{eq35}
	\overline{B}^{i}= \frac{B^{i}}{\sqrt{p}}\,,
\end{equation}\\[-12mm]
which gives\\[-8mm]
\begin{equation}\label{eq36}
	\overline{B}_{i}= \overline{g}_{ij}\overline{B}^{j}= \sqrt{p}B_{i}\,.
\end{equation}\\[-6mm]
Thus, we have
\begin{Theorem}
	Let $\overline{F}^{n}$ be the Finsler space obtained from $F^{n}$ by the \textsl{h-}Matsumoto change \eqref{eq3}. If $\overline{F}^{n-1}$ are the hypersurface of these spaces then the vector $b_{i}$ is tangential to the hypersurface $F^{n-1}$ if and only if every vector normal to $F^{n-1}$ is also normal to $\overline{F}^{n-1}$.
 \end{Theorem}
As $h_{ij}=g_{ij}-l_{i}l_{j}$, in view of equation \eqref{eq6}  and \eqref{eq31}, we get \\[-6mm]
\begin{equation}\label{eq37}
	h_{ij}B^{j}_{\alpha}B^{i}=0,\quad h_{ij}B^{i}=B_{j}.
\end{equation}
Then the tensors $B_{ij}$ and $Q_{i}$ given by \eqref{eq27}, satisfy the relations\\[-6mm]
\begin{equation}
		B_{ij}B^{i}B^{j}_{\alpha}=0, \quad B_{ij}B^{i}=B_{j}, \quad Q_{j}B^{j}=0.
\end{equation}
Transvecting \eqref{eq20} by $B_{i}$ and using $l^{i}B_{i}=0=m^{i}B_{i}$, we get\\[-6mm]
\begin{equation}
	\overline{g}^{is}B_{i}=qB^{s}.
\end{equation}\\[-10mm]
From \eqref{eq10}, \eqref{eq25}  and \eqref{eq36}, we get\\[-6mm]
\begin{equation*}
	\overline{H}_{\alpha}= \sqrt{p}(H_{\alpha}+ B_{i}D^{i}_{j}B^{j}_{\alpha}).
\end{equation*}\\[-9mm]
Contracting the above equation by $v^{\alpha}$ and using $v^{\alpha}B^{k}_{\alpha}=y^{k}$, we get\\[-6mm]
\begin{equation}\label{eq40}
	\overline{H}_{0}= \sqrt{p}(H_{0}+ B_{i}D^{i}).
\end{equation}\\[-9mm]
Equation \eqref{eq24} can be rewritten as \\[-4mm]
\begin{equation}\label{eq41}
D^{i}=\frac{1}{2} \Big\{\left[(q+q^{}_{1})p_{2}+q^{}_{2}p^{}_{3}m^{2}\right]E^{}_{00}+2q^{}_{2}p^{}_{2}LF_{\beta0}\Big\}l^{i}+\frac{1}{2} \Big\{\mu E_{00}+2q_{3}p_{2}LF_{\beta0}\Big\}m^{i}+ qp^{}_{2}LF^{i}_{0}.
\end{equation}
where\, $\mu= \left(qp^{}_{3}+q^{}_{2}p^{}_{2}+q^{}_{3}p^{}_{3}m^{2}\right)$ \,and \, $F_{\beta0}=F_{s0}m^{s}$.\\
Transvecting  the above equation  by $B_{i}$ and using $m_{i}B^{i}=0$  and  $l_{i}B^{i}=0$,  we get\\[-6mm]
\begin{equation}\label{eq42}
	D^{i}B_{i}= qp^{}_{2}LB_{i}F^{i}_{0}.
\end{equation}\\[-9mm]
Let the vector $b_{i}$ be gradiant, \textit{i.e.} $b_{i|j}=b_{j|i}$, then \\[-6mm]
\begin{equation}\label{eq43}
	F_{ij}=0.
\end{equation}\\[-9mm]
 M.$\,$K.$\,$Gupta and P.$\,$N.$\,$Pandey \cite{gupta2015finsler} proved the following Lemma,\\[-9mm] 
 \begin{Lemma}\label{L5}
 	If the \textsl{h-}vector $b_{i}$ is gradient then the scalar $\rho$ is constant.
 \end{Lemma}
From the above Lemma we get\\[-9mm]
\begin{equation}\label{eq44}
	\rho^{}_{i}=0.
\end{equation}\\[-9mm]
In view of \eqref{eq43}, the equation \eqref{eq42} becomes\\[-6mm]
 \begin{equation}\label{eq45}
 	D^{i}B_{i}=0. 
 \end{equation}\\[-9mm]
 and then equation \eqref{eq40} reduces to\\[-6mm]
\begin{equation}\label{eq46}
	\overline{H}_{0}= \sqrt{p}H_{0}.
\end{equation}\\[-9mm]
Thus in view of Lemma \eqref{L1}, we have
\begin{Theorem}
For the \textsl{h-}Matsumoto change, let the \textsl{h-}vector $b_{i}(x,y)$ be  gradient and tangent to the hypersurface $F^{n-1}$. Then the hypersurface $F^{n-1}$ is a hyperplane of the first kind if and only if hypersurface $\overline{F}^{n-1}$ is a hyperplane of the first kind.
\end{Theorem}
The second fundamental \textsl{h-}tensor $\overline{H}_{\alpha\beta}$ for hyperplane $\overline{F}^{n-1}$ is given by\\[-4mm]
\begin{equation*}
	\overline{H}_{\alpha\beta}=\overline{B}_{i}(B^{i}_{\alpha\beta}+\overline{F}^{i}_{jk}B^{j}_{\alpha}B^{k}_{\beta})+ \overline{M}_{\alpha}\overline{H}_{\beta}.
\end{equation*}
In view of equation \eqref{eq21}  and \eqref{eq36},  above equation gives\\[-4mm]
\begin{equation}\label{eq47}
 \overline{H}_{\alpha\beta}-\overline{M}_{\alpha}\overline{H}_{\beta}=\sqrt{p}(H_{\alpha\beta}+D^{i}_{jk}B_{i}B^{j}_{\alpha}B^{k}_{\beta})-\sqrt{p}{M}_{\alpha}{H}_{\beta}.
\end{equation}
 Using \eqref{eq43} and \eqref{eq44}, the equation \eqref{eq26} reduces to\\[-1mm]
\begin{equation}\begin{split}\label{eq48}
			D^{i}_{jk}&=\overline{g}^{is}\Big\{Q_{s}E_{kj}+pC_{jkm}D^{m}_{s}+V_{jkm}D^{m}_{s}-pC_{skm}D^{m}_{j}-V_{skm}D^{m}_{j}\\&\qquad \quad \quad -pC_{jsm}D^{m}_{k}-V_{jsm}D^{m}_{k}+B_{js}\beta_{k}+B_{sk}\beta_{j}-B_{jk}\beta_{s}\Big\}.
\end{split}\end{equation} 
Transvecting equation \eqref{eq48} by $B_{i}B^{j}_{\alpha}B^{k}_{\beta}$ and using\, $\overline{g}^{ij}B_{j}=qB^{i},\,\, B^{s}Q_{s}=0\,,\,\, B_{sk}B^{s}B^{k}_{\beta}=0$\,, we get\\[-8mm]
\begin{equation}\begin{split}\label{eq50}
	D^{i}_{jk}B_{i}B^{j}_{\alpha}B^{k}_{\beta}&=q{B}^{s}B^{j}_{\alpha}B^{k}_{\beta}\Big\{pC_{jkm}D^{m}_{s}+V_{jkm}D^{m}_{s}-pC_{skm}D^{m}_{j}-V_{skm}D^{m}_{j}\\&\qquad \qquad \qquad \qquad-pC_{jsm}D^{m}_{k}-V_{jsm}D^{m}_{k}-B_{jk}\beta_{s}\Big\}.
\end{split}\end{equation}
In view of equation \eqref{eq41}, and using the indicatory property of $C_{ijk}$, $V_{ijk}$, $h_{ij}$, $m_{i}$, the equation \eqref{eq25} can be rewritten as \\[-8mm]
\begin{equation}\label{eq52}
			D^{m}_{s}=\overline{g}^{mr} \Big\{\lambda h_{sr}+\phi m_{s}m_{r}+ Q_{r}E_{s0}\Big\},
\end{equation} \\[-10mm]
where\\[-6mm]
\begin{equation}\label{eq53}
		\lambda= \left[-\mu \left(\frac{p \rho}{L}+m^{2}K_{1}\right)+K_{1}\right]E_{00} \quad \text{and} \hspace{.9cm} \phi=\left[-\mu \left(2K_{1}+K_{2}m^{2}\right)+K_{2}\right]E_{00}.
\end{equation}
Transvecting equation \eqref{eq52} by $C_{jkm}$ we get\\[-6mm]
\begin{equation}
		C_{jkm}D^{m}_{s}= C_{jkm}\Big\{qg^{mr}+q^{}_{1}l^{m}l^{r}+q^{}_{2}\left(l^{m}m^{r}+l^{r}m^{m}\right)+q^{}_{3}m^{m}m^{r}\Big\} \Big\{\lambda h_{sr}+\phi m_{s}m_{r}+ Q_{r}E_{s0}\Big\},
	\end{equation}
which can be simplified as\\[-6mm]
\begin{equation}\label{eq54}
	C_{jkm}D^{m}_{s}= q\lambda C_{jsk}+ \Big[\left(q+q^{}_{3}m^{2}\right)\phi+q^{}_{3}\lambda\Big]\frac{\rho}{L}h_{jk}m_{s}+\frac{\rho}{L}\mu h_{jk}E_{s0}.
\end{equation}
Similarly we can write the expressions for $C_{skm}D^{m}_{j}$ and $C_{jsm}D^{m}_{k}$ as\\[-3mm]
\begin{equation}\label{eq55}
		C_{skm}D^{m}_{j}= q\lambda C_{jsk}+ \Big[\left(q+q^{}_{3}m^{2}\right)\phi+q^{}_{3}\lambda\Big]\frac{\rho}{L}h_{sk}m_{j}+\frac{\rho}{L}\mu h_{sk}E_{j0},
	\end{equation}\\[-8mm]
and\\[-8mm]
\begin{equation}\label{eq56}
	C_{jsm}D^{m}_{k}= q\lambda C_{jsk}+ \Big[\left(q+q^{}_{3}m^{2}\right)\phi+q^{}_{3}\lambda\Big]\frac{\rho}{L}h_{sj}m_{k}+\frac{\rho}{L}\mu h_{sj}E_{k0}.
\end{equation}
Transvecting equations \eqref{eq54}, \eqref{eq55} and \eqref{eq56} by $B^{s}B^{j}_{\alpha}B^{k}_{\beta}$, and using \eqref{eq12} and \eqref{eq37}, we get respectively\\[-8mm]
\begin{equation}\label{eq57}
		B^{s}C_{jkm}D^{m}_{s}B^{j}_{\alpha}B^{k}_{\beta}=q\lambda M_{\alpha \beta} + \frac{\rho}{L}\mu h_{jk}B^{j}_{\alpha}B^{k}_{\beta}B^{s}E_{s0}\,,
\end{equation}
\begin{equation}\label{eq58}
		B^{s}C_{skm}D^{m}_{j}B^{j}_{\alpha}B^{k}_{\beta}=q\lambda M_{\alpha \beta}\,,~~~~~~~~~~~~~~~~~~~~~~~~~~~~
\end{equation}
\begin{equation}\label{eq59}
		B^{s}C_{jsm}D^{m}_{k}B^{j}_{\alpha}B^{k}_{\beta}=q\lambda M_{\alpha \beta}\,.~~~~~~~~~~~~~~~~~~~~~~~~~~~~
\end{equation}\\[-8mm]
Again, transvecting Equation \eqref{eq52} by $V_{jkm}$, we get\\[-2mm]
\begin{equation}
		V_{jkm}D^{m}_{s}= V_{jkm}\Big\{qg^{mr}+\left(q^{}_{2}l^{r}+q^{}_{3}m^{r}\right)m^{m}\Big\}\Big\{\lambda h_{sr}+\phi m_{s}m_{r}+ Q_{r}E_{s0}\Big\}.
\end{equation}\\[-6mm]
By using Remark \ref{R3.1}, the above equation can be rewritten as \\[-2mm]
\begin{equation*}\begin{split}
		V_{jkm}D^{m}_{s}&=\Big[q\left\{K_{1}[h_{jk}m^{r}+h^{r}_{j}m_{k}+h^{r}_{k}m_{j}]+K_{2}m_{j}m_{k}m^{r}\right\}\\&\qquad\,\,+\left(q^{}_{2}l^{r}+q^{}_{3}m^{r}\right)\left[(2K_{1}+K_{2}m^{2})m_{j}m_{k}+K_{1}m^{2}h_{jk}\right]\Big]\Big\{\lambda h_{sr}+\phi m_{s}m_{r}+ Q_{r}E_{s0}\Big\},
	\end{split}\end{equation*}
which can be simplified as
\begin{equation}\begin{split}\label{eq60}
		V_{jkm}D^{m}_{s}&=\Big\{\psi_{1} K_{1}h_{jk}+(\psi_{1}K_{2}+2\psi_{2}K_{1})m_{j}m_{k}\Big\}m_{s}+q\lambda K_{1}\left(h_{js}m_{k}+h_{sk}m_{j}\right)\\& \qquad + \Big\{\mu \left[m^{2}K_{1}h_{jk}+(2K_{1}+K_{2}m^{2})m_{j}m_{k}\right]\Big\}E_{s0},
\end{split}\end{equation}
where\\[-8mm]
\begin{equation*}
	\psi^{}_{1}= \left(\lambda+\phi m^{2}\right) \left(q+q^{}_{3}m^{2}\right)\quad \text{and}\qquad \psi^{}_{2}= \left(\lambda+\phi m^{2}\right)q^{}_{3}+ q\phi.
\end{equation*}
Similarly we can write the expression for  $V_{skm}D^{m}_{j}$ and $V_{jsm}D^{m}_{k}$ as\\[-3mm]
\begin{equation}\begin{split}\label{eq61}
		V_{skm}D^{m}_{j}&= \Big\{\psi^{}_{1} K_{1}h_{sk}+(\psi^{}_{1}K_{2}+2\psi^{}_{2}K_{1})m_{s}m_{k}\Big\}m_{j}+q\lambda K_{1}\left(h_{js}m_{k}+h_{jk}m_{s}\right)\\& \qquad+ \Big\{\mu \left[m^{2}K_{1}h_{sk}+(2K_{1}+K_{2}m^{2})m_{s}m_{k}\right]\Big\}E_{j0},
\end{split}\end{equation}\\[-6mm]
and\\[-6mm]
\begin{equation}\begin{split}\label{eq62}
		V_{jsm}D^{m}_{k}&= \Big\{\psi_{1} K_{1}h_{sj}+(\psi^{}_{1}K_{2}+2\psi^{}_{2}K_{1})m_{j}m_{s}\Big\}m_{k}+q\lambda K_{1}\left(h_{sk}m_{j}+h_{jk}m_{s}\right)\\& \qquad+ \Big\{\mu \left[m^{2}K_{1}h_{sj}+(2K_{1}+K_{2}m^{2})m_{s}m_{j}\right]\Big\}E_{k0}.
\end{split}\end{equation}
Contracting equation \eqref{eq60}, \eqref{eq61}, \eqref{eq62} by $B^{s}B^{j}_{\alpha}B^{k}_{\beta}$\, and using $B^{i}m_{i}=0=h_{ij}B^{i}B^{j}_{\alpha}$\,, we get respectively\\[-8mm]
\begin{equation}\label{eq63}
		 ~~~~~~~~~~~B^{s}V_{jkm}D^{m}_{s}B^{j}_{\alpha}B^{k}_{\beta}\,=\,\mu \left\{K_{1}m^{2}h_{jk}+(2K_{1}+K_{2}m^{2})m_{j}m_{k}\right\}B^{j}_{\alpha}B^{k}_{\beta}B^{s}E_{s0},
\end{equation}
\begin{equation}\label{eq64}
B^{s}V_{skm}D^{m}_{j}B^{j}_{\alpha}B^{k}_{\beta}=0,~~~~~~~~~~~~~~~~~~~~~~~~~~~~~~~~~~~~~~~~~~~~~~~~~~~~~~~~~~
\end{equation}
\begin{equation}\label{eq65}
B^{s}V_{jsm}D^{m}_{k}B^{j}_{\alpha}B^{k}_{\beta}=0.~~~~~~~~~~~~~~~~~~~~~~~~~~~~~~~~~~~~~~~~~~~~~~~~~~~~~~~~~~
\end{equation}
putting the value of equation \eqref{eq57}, \eqref{eq58}, \eqref{eq59}, \eqref{eq63}, \eqref{eq64} and \eqref{eq65} in equation \eqref{eq50}, we get\\[-8mm]
\begin{equation}\label{eq66}
		D^{i}_{jk}B_{i}B^{j}_{\alpha}B^{k}_{\beta}= \Big[\frac{\mu \rho}{L}-q\left\{\left[(2K_{1}+K_{2}m^{2})m_{j}m_{k}\right]+K_{1}m^{2}h_{jk}\right\}+B_{jk}\Big]B^{j}_{\alpha}B^{k}_{\beta}B^{s}E_{s0}- q\lambda M_{\alpha \beta}.~~
\end{equation}
 Now taking the relative \textsl{h}-covariant differentiation of $b_{i}B^{i}=0$  with respect to the Cartan connection of $F^{n-1}$, we get\\[-12mm]
\begin{equation*}
	b_{i|j}B^{i}+b_{i}B^{i}_{|\beta}=0\,.
\end{equation*}\\[-9mm]
In view of equation \eqref{eq13} and \eqref{eq14}, above equation becomes\\[-6mm]
\begin{equation*}
	\Big(b_{i|j}B^{j}H_{\beta}+b_{i|j}B^{j}_{\beta}\Big)B^{i}- b_{i}H_{\alpha \beta}B^{\alpha}_{j}g^{ij}=0,
\end{equation*}\\[-6mm]
which on contraction by $v^{\beta}$ and using \eqref{eq11}, gives\\[-6mm]
\begin{equation*}
	b_{i|0}B^{i}= (H_{\alpha}+M_{\alpha}H_{0})B^{\alpha}_{j}b^{j}- b_{i|j}H_{0}B^{i}B^{j}.
\end{equation*}\\[-8mm]
In view of Lemma \ref{L1}, if the hypersurface to be first kind then $H_{0}=0=H_{\alpha}$. Thus the above equation reduces to $b_{i|0}B^{i}=0$. 
The vector $b_{i}$ is gradient, \textit{i.e.} $b_{i|j}=b_{j|i}$, then we get\\[-4mm]
\begin{equation}\label{eq67}
	E_{i0}B^{i}= b_{i|0}B^{i}=0.
\end{equation}
Therefore equation \eqref{eq66} reduces to\\[-6mm]
\begin{equation}\label{eq68}
	D^{i}_{jk}B_{i}B^{j}_{\alpha}B^{k}_{\beta}= -q\lambda M_{\alpha \beta}.
\end{equation}\\[-8mm]
In view of the above equation and \eqref{eq47}, we get\\[-6mm]
\begin{equation}\label{eq69}
	\overline{H}_{\alpha\beta}-\overline{M}_{\alpha}\overline{H}_{\beta}=\sqrt{p}(H_{\alpha\beta}-q\lambda M_{\alpha \beta})-\sqrt{p}{M}_{\alpha}{H}_{\beta}.
\end{equation}\\[-9mm]
Now transvecting \eqref{eq18}  by $B^{i}_{\alpha}B^{j}_{\beta}B^{k}$ and in view of equation \eqref{eq34}  and \eqref{eq37}, we obtain\\[-6mm]
\begin{equation}\label{eq38}
	\overline{C}_{ijk}B^{i}_{\alpha}B^{j}_{\beta}B^{k}= pC_{ijk}B^{i}_{\alpha}B^{j}_{\beta}B^{k}.
\end{equation}\\[-9mm]
From \eqref{eq12} and \eqref{eq35}, equation \eqref{eq38}  may be written as\\[-6mm]
\begin{equation}\label{eq70}
	\overline{M}_{\alpha\beta}=\sqrt{p}M_{\alpha\beta}.
\end{equation}\\[-9mm]
\noindent Thus from \eqref{eq69} and \eqref{eq70}, we have
\begin{Theorem} For the \textsl{h-}Matsumoto change, let the \textsl{h}-vector $b_{i}$ be a gradient and tangential to hypersurface $F^{n-1}$ and satisfies condition \eqref{eq67}. Then\\[-8mm]
	\begin{enumerate}
		\item {$\overline{F}^{n-1}$ is a hyperplane of second kind if $F^{n-1}$ is hyperplane of second kind and $M_{\alpha \beta}=0.$}
		\item{ $\overline{F}^{n-1}$} is a hyperplane of third kind if $F^{n-1}$ is hyperplane of third kind.
	\end{enumerate} 
\end{Theorem}
For the \textsl{h-}Matsumoto change, let the vector $b_{i}$ be parallel with respect to the Cartan connection of $F^{n}$. Then Cartan connection coefficient and Berwald connection coefficient for both the spaces are identical  \cite{first paper}, \textit{i.e.}\\[-8mm]
\begin{equation}\label{eq71}
	\overline{F}^{i}_{jk}= F^{i}_{jk},
\end{equation}\\[-9mm]
and\\[-8mm]
\begin{equation}\label{71}
	\overline{G}^{i}_{jk}= G^{i}_{jk}.
\end{equation}\\[-9mm]
\noindent Thus we have
\begin{Theorem}
	For the \textsl{h-}Matsumoto change, let the vector $b_{i}(x,y)$ be parallel with respect to the Cartan connection of $F^{n}$ and tangent to the hypersurface $F^{n-1}$. Then $\overline{F}^{n-1}$ is a hyperplane of the second(third) kind if and only if $F^{n-1}$ is also a hyperplane of the second(third) kind. \\[-5mm]
\end{Theorem}
\noindent The (\textsl{v})\textsl{hv-}torsion tensor for the hyperplane of first kind is given by \cite{Intrinsic}\\[-3mm]
\begin{equation}\label{eq72}
	P^{\alpha}_{\beta\gamma}= B^{\alpha}_{i}K^{i}_{\beta \gamma},
\end{equation}\\[-8mm]
where $K^{i}_{\beta \gamma}=P^{i}_{jk}B^{j}_{\alpha}B^{k}_{\beta}$.
Now contracting equation \eqref{3} by $K^{j}_{\beta\gamma}$ and using the above equation, we get\\[-8mm]
\begin{equation}\label{eq73}
	K^{i}_{\beta\gamma}= B^{i}_{\delta}P^{\delta}_{\beta\gamma}+B^{i}B_{h}K^{h}_{\beta\gamma}.
\end{equation}
Since the \textsl{(v)hv-}torsion tensor $P^{i}_{jk}$ is given by\\[-6mm]
\begin{equation}
	P^{i}_{jk}=G^{i}_{jk}-F^{i}_{jk},
\end{equation}\\[-6mm]
In view of \eqref{eq71} and \eqref{71}, the above equation gives\\[-6mm]
\begin{equation}\label{11}
	\overline{P}^{i}_{jk}= P^{i}_{jk}.
\end{equation}
 Thus in view of \eqref{11} and $K^{i}_{\beta \gamma}=P^{i}_{jk}B^{j}_{\alpha}B^{k}_{\beta}$, we get \, $\overline{K}^{i}_{\beta\gamma}=K^{i}_{\beta\gamma}$.\\
 From the above relation and using equation \eqref{eq73} we get\\[-4mm]
\begin{equation}
\overline{P}^{\alpha}_{\beta\gamma}= \overline{B}^{\alpha}_{i}\left[P^{\alpha}_{\beta \gamma}B^{i}_{\alpha}+K^{j}_{\beta\gamma}B^{i}B_{j}\right].
\end{equation}\\[-8mm]
In view of equation \eqref{eq8} and \eqref{eq72}, the above equation gives us\\[-6mm]
\begin{equation}
	\overline{P}^{\alpha}_{\beta\gamma}= P^{\alpha}_{\beta\gamma}.
\end{equation}\\[-9mm]
 We know that if the \textsl{(v)hv-}torsion tensor $P^{i}_{jk}$ vanishes then a Finsler space $F^{n}$ is called Landsberg space. Thus we have 
\begin{Theorem}
For the \textsl{h-}Matsumoto, let the \textsl{h-}vector $b_{i}(x,y)$ be parallel with respect to the Cartan connection of $F^{n}$ and tangent to the hypersurface $F^{n-1}$. Then a hyperplane $F^{n-1}$ of first kind is Landsberg space if and only if the hyperplane $\overline{F}^{n-1}$ of first kind is a Landsberg space.
\end{Theorem}
For the hyperplane of first kind, the Berwald connection coefficients $G^{\alpha}_{\beta\gamma}$ are given by \cite{Intrinsic}\\[-6mm]
\begin{equation}\label{eq74}
	G^{\alpha}_{\beta\gamma}= B^{\alpha}_{i}A^{i}_{\beta\gamma},
\end{equation}\\[-9mm]
where $A^{i}_{\beta\gamma}=B^{i}_{\beta\gamma}+G^{i}_{jk}B^{j}_{\beta}B^{k}_{\gamma}$.
Now contracting equation \eqref{3} by $A^{j}_{\beta\gamma}$ and using the above equation, we get\\[-9mm]
\begin{equation}\label{eq75}
	A^{i}_{\beta\gamma}=B^{i}_{\delta}G^{\delta}_{\beta\gamma}+B^{i}B_{h}A^{h}_{\beta\gamma}.
\end{equation}
 Thus in view of equation \eqref{eq71}, \eqref{71}, and from the above equation we get\, $\overline{A}^{i}_{\beta\gamma}={A}^{i}_{\beta\gamma}$.\\[1mm]
By using this relation and equation \eqref{eq74} and \eqref{eq75}, we obtain\\[-4mm]
\begin{equation}
	\overline{G}^{\alpha}_{\beta\gamma}= \overline{B}^{\alpha}_{i}\left(B^{i}_{\delta}G^{\delta}_{\beta\gamma}+B^{i}B_{h}A^{h}_{\beta\gamma}\right),
\end{equation}\\[-9mm]
In view of \eqref{eq8}, the above equation gives as\\[-6mm]
\begin{equation}\label{eq76}
	\overline{G}^{\alpha}_{\beta\gamma}={G}^{\alpha}_{\beta\gamma}.
\end{equation}
If the Berwald connection coefficients $G^{i}_{jk}$ are function of position only then a Finsler space $F^{n}$ is called Berwald space. Thus from \eqref{eq76}, we obtain
\begin{Theorem}
	For the \textsl{h-}Matsumoto change, let the \textsl{h-}vector $b_{i}(x,y)$ be parallel with respect to the Cartan connection of $F^{n}$ and tangent to the hypersurface $F^{n-1}$. Then a hyperplane $F^{n-1}$ of first kind is Berwald space if and only if the hyperplane $\overline{F}^{n-1}$ of first kind is a Berwald space.
\end{Theorem}
\section*{Discussion}
\textbf{`` $\overline{F}^{n-1}$ is a hyperplane of third kind if $F^{n-1}$ is hyperplane of third kind."}\\ 
\doublespacing
This result has been proved by Gupta and Gupta \cite{gupta2017h} for \textsl{h-}exponential change (which is infinite in nature) with the \textsl{h-}vector $b_{i}$ be gradient and tangential to hypersurface $\overline{F}^{n-1}$ and satisfies the condition \\[-14mm]
\begin{equation}\label{10}
	\beta_{r}C^{r}_{ij}=0.
\end{equation}\\[-12mm]
While, Gupta and Pandey \cite{Gupta h2018} have also obtained the same result with same condition for Kropina change (which is finite in nature) with an \textsl{h-}vector.\\
In the present paper we have proved the same result for \textsl{h-}Matsumoto change (which is also infinite  in nature) without using the condition \eqref{10}.\\
Notice that the above result holds\\[-14mm]
\begin{enumerate}
	\item[(i)]  For both the changes \textit{i.e.} \textsl{h-}exponential change\cite{gupta2017h} (infinite nature) and Kropina change\cite{Gupta h2018} (finite nature) with \textsl{h-}vector by using condition \eqref{10}.\\[-12mm]
	\item[(ii)] For \textsl{h-}Matsumoto change (infinite nature) without using \eqref{10}.\\[-12mm]
\end{enumerate}
  \textbf{The question is that}, Is there any specific change with \textsl{h-}vector (without using the  condition \eqref{10}) for which $\overline{F}^{n-1}$ is a hyperplane of third kind if $F^{n-1}$ is hyperplane of third kind?
\normalsize
\begin{spacing}{1.1}

\end{spacing}
\end{document}